\documentclass[12pt]{amsart}
\usepackage{eepic}
\def\N{{\Bbb N}}
\def\Z{{\Bbb Z}}

\newtheorem{Theorem}{Theorem}[section]
\newtheorem{Corollary}[Theorem]{Corollary}
\newtheorem{Lemma}[Theorem]{Lemma}

\newtheorem{Thm}{Theorem}
\newtheorem{Cor}[Thm]{Corollary}
\theoremstyle{definition}
\newtheorem{Definition}[Theorem]{Definition}
\newtheorem{Example}{Example}
\newtheorem*{Problem}{Problem}
\theoremstyle{remark}
\newtheorem*{Remark}{Remark}
\begin{document}
\sloppy
\title{
Coxeter systems with two-dimensional \\
Davis-Vinberg complexes
}
\author{Tetsuya Hosaka} 
\address{Department of Mathematics, Utsunomiya University, 
Utsunomiya, 321-8505, Japan}
\date{September 15, 2003}
\email{hosaka@cc.utsunomiya-u.ac.jp}
\maketitle
\begin{abstract}
In this paper, 
we study Coxeter systems with two-dimensional Davis-Vinberg complexes. 
We show that for a Coxeter group $W$, 
if $(W,S)$ and $(W,S')$ are Coxeter systems 
with two-dimensional Davis-Vinberg complexes, then 
there exists $S''\subset W$ such that 
$(W,S'')$ is a Coxeter system which is isomorphic to $(W,S)$ 
and the sets of reflections in $(W,S'')$ and $(W,S')$ coincide.
Hence 
the Coxeter diagrams of $(W,S)$ and $(W,S')$ have 
the same number of vertices, the same number of edges 
and the same multiset of edge-labels.
This is an extension of results of A.Kaul and 
N.Brady, J.P.McCammond, B.M\"uhlherr and W.D.Neumann. 
\end{abstract}

\section{Introduction and preliminaries}

The purpose of this paper is to study
Coxeter systems  with two-dimensional Davis-Vinberg complexes.
A {\it Coxeter group} is a group $W$ having a presentation
$$\langle \,S \, | \, (st)^{m(s,t)}=1 \text{ for } s,t \in S \,
\rangle,$$ 
where $S$ is a finite set and 
$m:S \times S \rightarrow \N \cup \{\infty\}$ is a function 
satisfying the following conditions:
\begin{enumerate}
\item[(1)] $m(s,t)=m(t,s)$ for each $s,t \in S$,
\item[(2)] $m(s,s)=1$ for each $s \in S$, and
\item[(3)] $m(s,t) \ge 2$ for each $s,t \in S$
such that $s\neq t$.
\end{enumerate}
The pair $(W,S)$ is called a {\it Coxeter system}.
For a Coxeter group $W$, a generating set $S'$ of $W$ is called 
a {\it Coxeter generating set for $W$} if $(W,S')$ is a Coxeter system.
In a Coxeter system $(W,S)$, 
the conjugates of elements of $S$ are called {\it reflections}. 
We note that 
the reflections depend on the Coxeter generating set $S$ and not just on 
the Coxeter group $W$. 
Let $(W,S)$ be a Coxeter system.
For a subset $T \subset S$, 
$W_T$ is defined as the subgroup of $W$ generated by $T$, 
and called a {\it parabolic subgroup}.
If $T$ is the empty set, then $W_T$ is the trivial group.

A {\it diagram} is an undirected graph $\Gamma$ 
without loops or multiple edges 
with a map $\text{Edges}(\Gamma)\rightarrow\{2,3,4,\ldots\}$ 
which assigns an integer greater than $1$ to each of its edges. 
Since such diagrams are used to define Coxeter systems, 
they are called {\it Coxeter diagrams}. 

Let $(W,S)$ and $(W',S')$ be Coxeter systems. 
Two Coxeter systems $(W,S)$ and $(W',S')$ are 
said to be {\it isomorphic}, 
if there exists a bijection 
$\psi:S\rightarrow S'$ such that 
$$m(s,t)=m'(\psi(s),\psi(t))$$ 
for each $s,t \in S$, where 
$m(s,t)$ and $m'(s',t')$ are the orders of $st$ in $W$ 
and $s't'$ in $W'$, respectively.

In general, a Coxeter group does not always determine 
its Coxeter system up to isomorphism.
Indeed some counter-examples are known.

\begin{Example}[{\cite[p.38 Exercise~8]{Bo}}, \cite{BMMN}]\label{ex1}
It is known that 
the Coxeter groups defined 
by the diagrams in Figure~\ref{fig1} are 
isomorphic and $D_6$.
\begin{figure}[htbp]
\unitlength = 0.9mm
\begin{center}
\begin{picture}(80,28)(-40,-5)
\put(-25,0){\line(-1,0){20}}
\put(25,0){\line(1,0){20}}
\put(25,0){\line(10,17){10}}
\put(45,0){\line(-10,17){10}}
\put(-25,0){\circle*{1.3}}
\put(-45,0){\circle*{1.3}}
\put(25,0){\circle*{1.3}}
\put(45,0){\circle*{1.3}}
\put(35,17){\circle*{1.3}}
{\small
\put(-36,-4){$6$}
\put(34,-4){$3$}
\put(27,9){$2$}
\put(41,9){$2$}
}
\end{picture}
\end{center}
\caption[Two distinct Coxeter diagrams for $D_6$]{Two distinct Coxeter diagrams for $D_6$}
\label{fig1}
\end{figure}
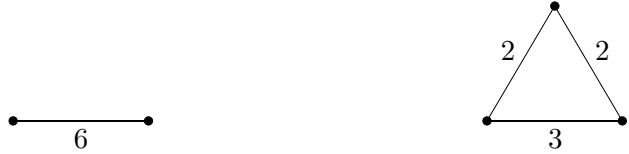
\end{Example}

\begin{Example}[\cite{Mu}, \cite{BMMN}]\label{ex2}
In \cite{Mu}, B.M\"uhlherr showed that
the Coxeter groups defined by 
the diagrams in Figure~\ref{fig2} are isomorphic.
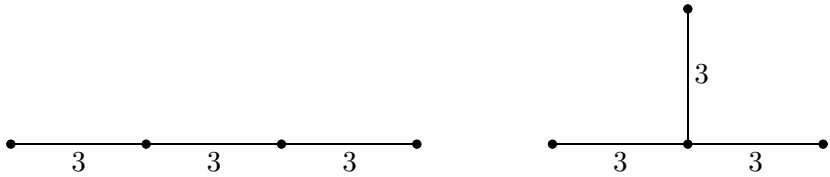
\begin{figure}[htbp]
\unitlength = 0.9mm
\begin{center}
\begin{picture}(140,30)(-70,-5)
\put(-5,0){\line(-1,0){60}}
\put(15,0){\line(1,0){40}}
\put(35,0){\line(0,1){20}}
\put(-5,0){\circle*{1.3}}
\put(-25,0){\circle*{1.3}}
\put(-45,0){\circle*{1.3}}
\put(-65,0){\circle*{1.3}}
\put(15,0){\circle*{1.3}}
\put(35,0){\circle*{1.3}}
\put(55,0){\circle*{1.3}}
\put(35,20){\circle*{1.3}}
{\small
\put(-16,-4){$3$}
\put(-36,-4){$3$}
\put(-56,-4){$3$}
\put(24,-4){$3$}
\put(44,-4){$3$}
\put(36,9){$3$}
}
\end{picture}
\end{center}
\caption[Coxeter diagrams for isomorphic Coxeter groups]{Coxeter diagrams for isomorphic Coxeter groups}
\label{fig2}
\end{figure}
\end{Example}

Here there exists the following natural problem.

\begin{Problem}[\cite{BMMN}, \cite{CD}]
When does a Coxeter group determine its Coxeter system up to isomorphism?
\end{Problem}

Recently, 
B.M\"uhlherr and R.Weidmann proved that 
skew-angled Coxeter systems are reflection rigid up to 
diagram twisting (\cite{MW}).

It is known that each Coxeter system 
$(W,S)$ defines a CAT(0) geodesic space $\Sigma (W,S)$ 
called the Davis-Vinberg complex (\cite{D1}, \cite{D3}, \cite{D2}, \cite{M}). 
Here $\dim \Sigma (W,S) \ge 1$ by definition, 
and $\dim \Sigma (W,S)=1$ if and only if 
the Coxeter group $W$ is isomorphic to the free product of some $\Z_2$.
Hence if $\dim \Sigma (W,S)=1$, then 
the Coxeter group $W$ is {\it rigid}, i.e., 
$W$ determines its Coxeter system up to isomorphism.
In this paper, 
we investigate Coxeter systems 
with two-dimensional Davis-Vinberg complexes.

\begin{Remark}
Let $(W,S)$ be a Coxeter system.
We note that 
$\dim \Sigma(W,S) \le 2$ if and only if 
$W_T$ is infinite for each $T\subset S$ such that $|T|>2$.
It is known that for $\{s_1,s_2,s_3\}\subset S$ if 
\begin{enumerate}
\item[(1)] $m(s_i,s_j)\ge 3$ for each $i,j\in\{1,2,3\}$ such that $i\neq j$ or
\item[(2)] $m(s_i,s_j)=\infty$ for some $i,j\in\{1,2,3\}$,
\end{enumerate}
then the parabolic subgroup 
$W_{\{s_1,s_2,s_3\}}$ is infinite (see \cite{Bo}).
Hence, for example, if 
\begin{enumerate}
\item[(1)] $(W,S)$ is of type $K_n$ (cf.\ \cite{K}), 
\item[(2)] all edge-labels of the Coxeter diagram of $(W,S)$ are odd, 
\item[(3)] all edge-labels of the Coxeter diagram of $(W,S)$ are 
greater than $2$ (i.e.\ $(W,S)$ is skew-angled), or 
\item[(4)] the Coxeter diagram of $(W,S)$ is tree, 
\end{enumerate}
then $\dim \Sigma (W,S) \le 2$.
\end{Remark}

We first recall some basic properties of Coxeter groups 
and Davis-Vinberg complexes in Section~2.
After some preliminaries in Section~3, 
we prove the following theorem in Section~4.

\begin{Thm}\label{thm1}
Let $(W,S)$ and $(W',S')$ be Coxeter systems 
with two-dimensional Davis-Vinberg complexes. 
Suppose that there exists an isomorphism $\phi:W\rightarrow W'$.
For each $s\in S$, 
if $\phi(s)$ is not a reflection in $(W',S')$, then 
there exist unique elements $t\in S$ and $s',t'\in S'$ 
such that for some $w'\in W'$, 
\begin{enumerate}
\item[(1)] $m(s,t)=2$, 
\item[(2)] $m(s,u)=\infty$ for each $u\in S\setminus\{s,t\}$, 
\item[(3)] $\phi(W_{\{s,t\}})=w'W'_{\{s',t'\}}(w')^{-1}$, 
\item[(4)] $m'(s',t')=2$, 
\item[(5)] $m'(s',u')=\infty$ for each $u'\in S'\setminus\{s',t'\}$, 
\item[(6)] $\phi(s)=w's't'(w')^{-1}$ and
\item[(7)] $\phi(t)=w't'(w')^{-1}$.
\end{enumerate}
\end{Thm}

Here 
we can define an automorphism $\psi$ of $W$ as follows:
for each $s\in S$, 
\begin{enumerate}
\item[(1)] if $\phi(s)$ is a reflection in $(W',S')$, then $\psi(s)=s$, and 
\item[(2)] if $\phi(s)$ is not a reflection in $(W',S')$, then $\psi(s)=st$, 
where $t$ is a unique element of $S$ such that $m(s,t)=2$.
\end{enumerate}
Then the Coxeter systems $(W,S)$ and $(W,\psi(S))$ are isomorphic 
and by Theorem~\ref{thm1} 
the isomorphism $\phi:W\rightarrow W'$ maps reflections in $(W,\psi(S))$ 
onto reflections in $(W',S')$.
Thus we obtain the following theorem.

\begin{Thm}\label{thm2}
Let $(W,S)$ and $(W,S')$ be Coxeter systems 
with two-dimensional Davis-Vinberg complexes. 
Then there exists $S''\subset W$ such that 
$(W,S'')$ is a Coxeter system which is isomorphic to $(W,S)$ 
and the sets of reflections in $(W,S'')$ and $(W,S')$ coincide.
\end{Thm}

This implies the following corollary which is an extension of results of 
A.Kaul (\cite{K}) and 
N.Brady, J.P.McCammond, B.M\"uhlherr and W.D.Neumann 
(\cite[Lemma~5.3]{BMMN}).

\begin{Cor}\label{cor2}
For a Coxeter group $W$, 
if $(W,S)$ and $(W,S')$ are Coxeter systems 
with two-dimensional Davis-Vinberg complexes, then 
the Coxeter diagrams of $(W,S)$ and $(W,S')$ have 
the same number of vertices, the same number of edges 
and the same multiset of edge-labels.
\end{Cor}

Here a {\it multiset} is a collection in which the order of 
the entries does not matter, but multiplicities do.
Thus the multisets $\{1,1,2\}$ and $\{1,2,2\}$ are different.
In Corollary~\ref{cor2}, we can not omit the assumption 
``with two-dimensional Davis-Vinberg complexes''
by Example~\ref{ex1}.

\section{Basics on Coxeter groups and Davis-Vinberg complexes}

In this section, we introduce some basic properties 
of Coxeter groups and Davis-Vinberg complexes.

\begin{Definition}
Let $(W,S)$ be a Coxeter system and $T \subset S$.
The subset $T$ is called a {\it spherical subset of $S$}, 
if the parabolic subgroup $W_T$ is finite.
\end{Definition}

\begin{Definition}
Let $(W,S)$ be a Coxeter system and $w\in W$.
A representation $w=s_1\cdots s_l$ ($s_i \in S$) is said to be 
{\it reduced}, if $\ell(w)=l$, 
where $\ell(w)$ is the minimum length of 
word in $S$ which represents $w$.
\end{Definition}

The following lemmas are known.

\begin{Lemma}[\cite{Bo}, \cite{Br}, \cite{D1}, \cite{Hu}]\label{lem1}
Let $(W,S)$ be a Coxeter system.
\begin{enumerate}
\item[(i)] Let $w\in W$ and 
let $w=s_1\cdots s_l$ be a representation.
If $\ell(w)<l$, then $w=s_1\cdots \hat{s_i} \cdots \hat{s_j} \cdots s_l$ 
for some $1\le i<j\le l$. 
\item[(ii)] Let $w\in W$ and 
let $w=s_1\cdots s_l$ be a representation.
Then the length $\ell(w)$ is even if and only if $l$ is even.
\item[(iii)] For each $w\in W$, 
there exists a unique subset $S(w)\subset S$ such that 
$S(w)=\{s_1,\dots,s_l\}$ for every reduced representation 
$w=s_1\cdots s_l$ ($s_i \in S$). 
\item[(iv)] Let $w\in W$ and $T \subset S$. 
Then $w\in W_T$ if and only if $S(w)\subset T$.
\item[(v)] For each subset $T \subset S$, 
$(W_T,T)$ is a Coxeter system. 
\item[(vi)] For each subsets $T_1, T_2 \subset S$, 
$W_{T_1}=W_{T_2}$ if and only if $T_1=T_2$.
\item[(vii)] If $W$ is finite, 
then there exists a unique element $w_0\in W$ of longest length. 
\end{enumerate}
\end{Lemma}

\begin{Lemma}[{\cite{Bo}, \cite[Lemma~7.11]{D1}}]\label{lem2}
Let $(W,S)$ be a Coxeter system, 
let $T\subset S$ and let $w\in W_T$. 
Then the following statements are equivalent:
\begin{enumerate}
\item[{\rm (1)}] $W_T$ is finite and 
$w$ is the element of longest length in $W_T$;
\item[{\rm (2)}] $\ell(wt)<\ell(w)$ for each $t\in T$.
\end{enumerate}
\end{Lemma}

\begin{Lemma}[{\cite[p.12 Proposition~3]{Bo}}]\label{lem2-5}
Let $(W,S)$ be a Coxeter system and let $s,t\in S$.
Then $s$ is conjugate to $t$ if and only if 
there exists a sequence $s_1,\dots,s_n\in S$ such that 
$s_1=s$, $s_n=t$ and $m(s_i,s_{i+1})$ is odd for each $i\in\{1,\dots,n-1\}$.
\end{Lemma}

\begin{Lemma}[{\cite[Results~3.7 and 3.8]{BMMN}}]\label{lem3-1}
Let $(W,S)$ and $(W,S')$ be Coxeter systems. Then 
\begin{enumerate}
\item[(1)] if $S\subset R_{S'}$ then $R_S=R_{S'}$, and
\item[(2)] if $R_S=R_{S'}$ then $|S|=|S'|$,
\end{enumerate}
where $R_S$ and $R_{S'}$ are the sets 
of all reflections in $(W,S)$ and $(W,S')$, respectively.
\end{Lemma}

By Results~1.8, 1.9 and 1.10 in \cite{BMMN}, 
we obtain the following theorem.

\begin{Theorem}[cf.\ \cite{BMMN}]\label{thm4-0}
Let $(W,S)$ and $(W',S')$ be Coxeter systems. 
Suppose that there exists an isomorphism $\phi:W\rightarrow W'$. 
Then for each maximal spherical subset $T\subset S$, 
there exists a unique maximal spherical subset $T'\subset S'$ 
such that $\phi(W_T)=w'W'_{T'}(w')^{-1}$ for some $w'\in W'$. 
\end{Theorem}

We introduce a definition of the Davis-Vinberg complex 
of a Coxeter system.

\begin{Definition}[\cite{D1}, \cite{D3}, \cite{D2}]
Let $(W,S)$ be a Coxeter system and 
let $W{\mathcal{S}}^f$ denote the set of all left cosets of the form
$w W_T$, with $w \in W$ and a spherical subset $T \subset S$.
The set $W{\mathcal{S}}^f$ is partially ordered by inclusion.
The Davis-Vinberg complex $\Sigma(W,S)$ is defined as 
the geometric realization of the partially ordered set $W{\mathcal{S}}^f$
(\cite{D1}, \cite{D3}).
Here it is known that 
$\Sigma(W,S)$ has a structure of 
a PE (i.e.\ Piecewise Euclidean) cell complex 
whose $1$-skeketon is the Cayley graph of $W$ with respect to $S$ (\cite{D2}).
Then the vertex set of 
each cell of the PE cell complex $\Sigma(W,S)$ is $wW_T$ 
for some $w \in W$ and some spherical subset $T$ of $S$.
The Coxeter group $W$ acts 
properly discontinuously and cocompactly as isometries 
on the PE cell complex $\Sigma(W,S)$ with the natural metric 
(\cite{D1}, \cite{D2}).
\end{Definition}

\begin{Remark}
For a Coxeter system $(W,S)$, 
by the definition of $\Sigma(W,S)$, 
$$ \dim \Sigma(W,S)=
\max \{ |T| : \text{$T$ is a spherical subset of $S$}\}.$$
\end{Remark}

Theorem~\ref{thm4-0} implies the following lemma.

\begin{Lemma}\label{lem4}
Let $(W,S)$ and $(W',S')$ be Coxeter systems 
with two-dimensional Davis-Vinberg complexes. 
Suppose that there exists an isomorphism $\phi:W\rightarrow W'$. 
\begin{enumerate}
\item[(i)]
For each two elements $s,t\in S$ such that $m(s,t)<\infty$, 
there exist unique two elements $s',t'\in S'$ 
such that $\phi(W_{\{s,t\}})=w'W'_{\{s',t'\}}(w')^{-1}$ 
(hence $m(s,t)=m'(s',t')$) for some $w'\in W'$.
\item[(ii)] 
The multisets of edge-labels of the Coxeter diagrams 
of $(W,S)$ and $(W',S')$ coincide.
\end{enumerate}
\end{Lemma}

\section{Lemmas on Coxeter groups}

We show some lemmas needed later.

\begin{Lemma}\label{lem3}
Let $(W,S)$ be a Coxeter system and let $w\in W$. 
Suppose that $w^2=1$ and 
$\ell(w)=\min\{\ell(vwv^{-1}): v\in W\}$. 
Then $W_{S(w)}$ is finite and 
$w$ is the element of longest length in $W_{S(w)}$, 
where $S(w)$ is the subset of $S$ defined in Lemma~\ref{lem1}~(iii).
\end{Lemma}

\begin{proof}
Let $w=s_1\cdots s_l$ be a reduced representation. 
Since $w^2=1$, 
$$ s_1\cdots s_l=w=w^{-1}=s_l\cdots s_1.$$
Hence $\ell(ws_1)<\ell(w)$. 
By Lemma~\ref{lem1}~(i), 
$$ ws_1=(s_1\cdots s_l)s_1=s_1\cdots \hat{s_i} \cdots s_l $$
for some $i\in\{1,\dots,l\}$.
Suppose that $1<i\le l$. Then 
$$ s_1 w s_1 = s_2\cdots \hat{s_i} \cdots s_l, $$
and $\ell(s_1ws_1)<\ell(w)$. 
This contradicts the assumption 
$$\ell(w)=\min\{\ell(vwv^{-1}) : v\in W\}. $$
Thus $i=1$ and $ws_1=s_2\cdots s_l$.
Hence $w=(s_2\cdots s_l)s_1$ is reduced.

By iterating the above argument, 
$$ w=(s_{i+1}\cdots s_l)(s_1\cdots s_i) $$
is reduced for each $i\in\{1,\dots,l-1\}$. 
Hence $\ell(ws_i)<\ell(w)$ for each $i\in\{1,\dots,l\}$, i.e., 
$\ell(ws)<\ell(w)$ for each $s\in S(w)$. 
By Lemma~\ref{lem2}, 
$W_{S(w)}$ is finite and 
$w$ is the element of longest length in $W_{S(w)}$. 
\end{proof}

\begin{Remark}
Let $(W,S)$ and $(W',S')$ be Coxeter systems 
with two-dimensional Davis-Vinberg complexes. 
Suppose that there exists an isomorphism $\phi:W\rightarrow W'$.
Let $s\in S$. 
Since $(\phi(s))^2=1$, by Lemma~\ref{lem3}, either 
\begin{enumerate}
\item[(1)] $\phi(s)$ is a reflection in $(W',S')$, or
\item[(2)] $\phi(s)=w'(s't')^{m'}(w')^{-1}$ 
for some $w'\in W'$ and $s',t'\in S'$, where 
$m'(s',t')$ is even and $m'=m'(s',t')/2$. 
\end{enumerate}
\end{Remark}

\begin{Lemma}\label{lem6}
Let $(W,S)$ be a Coxeter system 
with two-dimensional Davis-Vinberg complex, 
let $s,t,a,b\in S$ and let $w,x\in W$. 
Suppose that $m(s,t)$ is even,
$m(a,b)$ is finite and 
$w(st)^m w^{-1}\in xW_{\{a,b\}}x^{-1}$, 
where $m=m(s,t)/2$.
Then $wW_{\{s,t\}}w^{-1}=xW_{\{a,b\}}x^{-1}$ and $\{s,t\}=\{a,b\}$. 
\end{Lemma}

\begin{proof}
Suppose that $m(s,t)$ is even,
$m(a,b)$ is finite and 
$w(st)^m w^{-1}\in xW_{\{a,b\}}x^{-1}$, 
where $m=m(s,t)/2$.
Let $v=w(st)^m w^{-1}$ and 
let $C$ and $D$ be the 2-cells in $\Sigma(W,S)$ 
such that 
$$C^{(0)}=wW_{\{s,t\}} \text{ and } D^{(0)}=xW_{\{a,b\}}.$$
Then $v$ is an isometry of $\Sigma(W,S)$ and 
the barycenter of $C$ is the unique fixed point of $v$ 
because $m(s,t)=2m$ and $\dim\Sigma(W,S)=2$ (cf.\ Figure~\ref{fig3}).
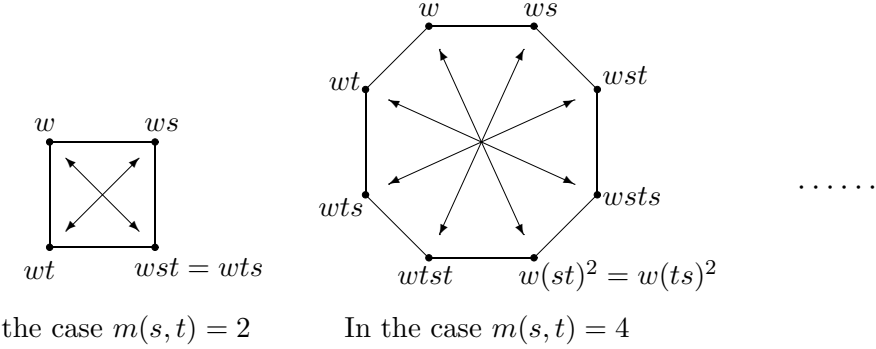
\begin{figure}[tbp]
\unitlength = 0.7mm
\begin{center}
\begin{picture}(180,60)(-100,-12)
\put(-10,0){\line(1,0){20}}
\put(10,0){\line(1,1){12}}
\put(22,12){\line(0,1){20}}
\put(22,32){\line(-1,1){12}}
\put(10,44){\line(-1,0){20}}
\put(-10,44){\line(-1,-1){12}}
\put(-22,32){\line(0,-1){20}}
\put(-22,12){\line(1,-1){12}}
\put(-82,2){\line(1,0){20}}
\put(-62,2){\line(0,1){20}}
\put(-62,22){\line(-1,0){20}}
\put(-82,22){\line(0,-1){20}}
\put(-10,0){\circle*{1.3}}
\put(10,0){\circle*{1.3}}
\put(22,12){\circle*{1.3}}
\put(22,32){\circle*{1.3}}
\put(10,44){\circle*{1.3}}
\put(-10,44){\circle*{1.3}}
\put(-22,32){\circle*{1.3}}
\put(-22,12){\circle*{1.3}}
\put(-82,2){\circle*{1.3}}
\put(-62,2){\circle*{1.3}}
\put(-62,22){\circle*{1.3}}
\put(-82,22){\circle*{1.3}}
\put(60,12){$\cdots\cdots$}
\put(0,22){\line(5,11){8}}
\put(0,22){\line(-5,-11){8}}
\put(0,22){\line(-5,11){8}}
\put(0,22){\line(5,-11){8}}
\put(0,22){\line(11,5){17.6}}
\put(0,22){\line(-11,-5){17.6}}
\put(0,22){\line(-11,5){17.6}}
\put(0,22){\line(11,-5){17.6}}
\put(8,39.6){\vector(1,2){0}}
\put(-8,4.4){\vector(-1,-2){0}}
\put(17.6,30){\vector(2,1){0}}
\put(-17.6,14){\vector(-2,-1){0}}
\put(17.6,14){\vector(2,-1){0}}
\put(-17.6,30){\vector(-2,1){0}}
\put(-8,39.6){\vector(-1,2){0}}
\put(8,4.4){\vector(1,-2){0}}
\put(-72,12){\vector(1,1){7}}
\put(-72,12){\vector(-1,-1){7}}
\put(-72,12){\vector(-1,1){7}}
\put(-72,12){\vector(1,-1){7}}
{\small
\put(-12,46){$w$}
\put(8,46){$ws$}
\put(23,33){$wst$}
\put(23,10){$wsts$}
\put(7,-5){$w(st)^2=w(ts)^2$}
\put(-16,-5){$wtst$}
\put(-31,8){$wts$}
\put(-29,32){$wt$}
\put(-85,24){$w$}
\put(-64,24){$ws$}
\put(-87,-4){$wt$}
\put(-66,-3.5){$wst=wts$}
\put(-98,-15){In the case $m(s,t)=2$}
\put(-26,-15){In the case $m(s,t)=4$}
}
\end{picture}
\end{center}
\caption[The isometry $v$ on the $2$-cell $C$]{The isometry $v$ on the $2$-cell $C$}
\label{fig3}
\end{figure}
Since 
$$v=w(st)^m w^{-1}\in xW_{\{a,b\}}x^{-1}, $$
there exists $u\in W_{\{a,b\}}$ such that $v=xux^{-1}$.
Then 
$$ v(xW_{\{a,b\}})=xux^{-1}(xW_{\{a,b\}})=x(uW_{\{a,b\}})=xW_{\{a,b\}}.$$
Hence $vD=D$.
In general, for each cell $E$ of $\Sigma(W,S)$ and each $y\in W$, 
if $yE=E$ then the isometry $y$ fixes the barycenter of $E$
by the definition of $\Sigma(W,S)$.
Thus the barycenter of $D$ is a fixed point of $v$.
On the other hand, 
the barycenter of $C$ is the unique fixed point of $v$. 
Hence $C=D$ and 
$$wW_{\{s,t\}}=C^{(0)}=D^{(0)}=xW_{\{a,b\}}.$$
Since $x^{-1}wW_{\{s,t\}}=W_{\{a,b\}}$, 
we have that 
$$x^{-1}w\in x^{-1}w W_{\{s,t\}}=W_{\{a,b\}}.$$
Hence $W_{\{s,t\}}=W_{\{a,b\}}$ and $\{s,t\}=\{a,b\}$
by Lemma~\ref{lem1}~(vi).
Since 
$$x^{-1}w, (x^{-1}w)^{-1}\in W_{\{s,t\}}=W_{\{a,b\}},$$ 
$x^{-1}w W_{\{s,t\}} w^{-1}x=W_{\{a,b\}}$. 
Hence we obtain that 
$wW_{\{s,t\}}w^{-1}=xW_{\{a,b\}}x^{-1}$. 
\end{proof}

\begin{Lemma}\label{lem5}
Let $(W,S)$ and $(W',S')$ be Coxeter systems 
with two-dimensional Davis-Vinberg complexes such that 
there exists an isomorphism $\phi:W\rightarrow W'$, 
let $s\in S$, let $s',t'\in S'$ and let $w'\in W'$.
Suppose that $m'(s',t')$ is even and 
$\phi(s)= w'(s't')^{m'}(w')^{-1}$, where $m'=m'(s',t')/2$.
Then there exists a unique element $t\in S$ such that 
\begin{enumerate}
\item[(1)] $\phi(W_{\{s,t\}})=w'W'_{\{s',t'\}}(w')^{-1}$, 
\item[(2)] $\phi(t)$ is a reflection in $(W',S')$, and 
\item[(3)] $m(s,t)=m'(s',t')=2$.
\end{enumerate}
\end{Lemma}

\begin{proof}
Suppose that $m'(s',t')$ is even and 
$\phi(s)= w'(s't')^{m'}(w')^{-1}$, where $m'=m'(s',t')/2$.
By Lemma~\ref{lem4}, 
there exist $r,t\in S$ and $x\in W$ such that 
$$\phi^{-1}(W'_{\{s',t'\}})=xW_{\{r,t\}}x^{-1}.$$
Here we note that $m(r,t)=m'(s',t')$. 

We first show that we may suppose $r=s$.

Since $\phi(s)=w'(s't')^{m'}(w')^{-1}$, 
$$ (\phi^{-1}(w'))^{-1} s \phi^{-1}(w')=\phi^{-1}((s't')^{m'})
\in \phi^{-1}(W'_{\{s',t'\}})=xW_{\{r,t\}}x^{-1}.$$ 
Hence $(\phi^{-1}(w'))^{-1} s \phi^{-1}(w')=xyx^{-1}$ 
for some $y\in W_{\{r,t\}}$.
Since the length of 
$(\phi^{-1}(w'))^{-1} s \phi^{-1}(w')$ is odd, 
the length of $y$ is also odd.
Hence $y$ is conjugate to either $r$ or $t$
because $y\in W_{\{r,t\}}$.
Here we may suppose that $y$ is conjugate to $r$. 
Then $s$ is conjugate to $r$, since 
$(\phi^{-1}(w'))^{-1} s \phi^{-1}(w')=xyx^{-1}$. 

Now we show that $s=r$. 
If $s\neq r$, then 
there exists $a\in S\setminus \{s\}$ 
such that $m(s,a)$ is odd by Lemma~\ref{lem2-5}.
By Lemma~\ref{lem4}, 
there exist $a',b'\in S'$ such that 
$\phi(W_{\{s,a\}})=x'W'_{\{a',b'\}}(x')^{-1}$ 
for some $x'\in W'$. 
Here we note that $m(s,a)=m'(a',b')$.
Then 
$$ w'(s't')^{m'}(w')^{-1}=\phi(s)\in 
\phi(W_{\{s,a\}})=x'W'_{\{a',b'\}}(x')^{-1}. $$
Hence $\{s',t'\}=\{a',b'\}$ by Lemma~\ref{lem6} 
and $m'(a',b')=m'(s',t')$.
Then 
$$m(s,a)=m'(a',b')=m'(s',t'). $$
Here $m(s,a)$ is odd.
This contradicts the assumption $m'(s',t')=2m'$ is even.
Thus $s=r$.

Then $\phi^{-1}(W'_{\{s',t'\}})=xW_{\{s,t\}}x^{-1}$ and 
$$ w'(s't')^{m'}(w')^{-1}=\phi(s)\in 
\phi(W_{\{s,t\}})=(\phi(x))^{-1}W'_{\{s',t'\}}\phi(x). $$
By Lemma~\ref{lem6}, 
$(\phi(x))^{-1}W'_{\{s',t'\}}\phi(x)=w'W'_{\{s',t'\}}(w')^{-1}$.
Hence 
$$\phi(W_{\{s,t\}})=w'W'_{\{s',t'\}}(w')^{-1}.$$
Here we note that such $t\in S$ is unique by Lemma~\ref{lem4}.

Next we show that $\phi(t)$ is a reflection.
Here $\phi(t)$ is a reflection if and only if 
the length $\ell(\phi(t))$ is odd, 
because 
$\phi(t)\in \phi(W_{\{s,t\}})=w'W'_{\{s',t'\}}(w')^{-1}$.
Now we suppose that the length $\ell(\phi(t))$ is even.
Then the lengthes of $\phi(s)= w'(s't')^{m'}(w')^{-1}$ 
and $\phi(t)$ are even and 
the set $\{\phi(s),\phi(t)\}$ generates 
$\phi(W_{\{s,t\}})=w'W'_{\{s',t'\}}(w')^{-1}$.
In general, for $f,g\in W$ 
if $\ell(f)$ and $\ell(g)$ are even, 
then the length $\ell(fg)$ is even
by Lemma~\ref{lem1}~(ii).
Hence the length of each element of 
$\phi(W_{\{s,t\}})$ is even.
On the other hand, 
the length of 
$w's'(w')^{-1}\in w'W'_{\{s',t'\}}(w')^{-1}=\phi(W_{\{s,t\}})$ is odd.
This is a contradiction.
Thus the length of $\phi(t)$ is odd 
and $\phi(t)$ is a reflection.

Since $\phi(t)$ is a reflection, 
$\phi(t)=w'(s't')^ks'(w')^{-1}$ for some $0\le k<2m'$.
Then 
\begin{align*}
\phi(s)\phi(t)&=(w'(s't')^{m'}(w')^{-1})(w'(s't')^ks'(w')^{-1}) \\
&=w'(s't')^{m'}(s't')^ks'(w')^{-1} \\
&= w'(s't')^{m'+k}s'(w')^{-1}. 
\end{align*}
Hence $\phi(s)\phi(t)$ is a reflection and 
$(\phi(s)\phi(t))^2=1$, i.e., $(st)^2=1$.
This means that $m(s,t)=m'(s',t')=2$.
\end{proof}

\begin{Lemma}\label{lem7}
Let $(W,S)$ be a Coxeter system and let $s,t\in S$.
Suppose that $m(s,t)=2$ and 
$m(s,u)=\infty$ for each $u\in S\setminus\{s,t\}$.
Let $S'=(S\setminus\{s\})\cup \{st\}$.
Then
$(W,S')$ is a Coxeter system which is isomorphic to $(W,S)$.
\end{Lemma}

\begin{proof}
The map $\psi:S\rightarrow S'$ defined by 
$\psi(s)=st$ and $\psi(u)=u$ for each $u\in S\setminus \{s\}$ 
induces an automorphism 
$\psi:W\rightarrow W$, 
and $(W,S)$ and $(W,S')$ are isomorphic.
\end{proof}

\section{Proof of the main results}

Using some lemmas in Sections~2 and 3, 
we prove the main results.

\begin{Theorem}\label{Thm1}
Let $(W,S)$ and $(W',S')$ be Coxeter systems 
with two-dimensional Davis-Vinberg complexes. 
Suppose that there exists an isomorphism $\phi:W\rightarrow W'$.
For each $s\in S$, 
if $\phi(s)$ is not a reflection in $(W',S')$, then 
there exist unique $t\in S$ and $s',t'\in S'$ 
such that for some $w'\in W'$, 
\begin{enumerate}
\item[(1)] $m(s,t)=2$, 
\item[(2)] $m(s,u)=\infty$ for each $u\in S\setminus\{s,t\}$, 
\item[(3)] $\phi(W_{\{s,t\}})=w'W'_{\{s',t'\}}(w')^{-1}$, 
\item[(4)] $m'(s',t')=2$, 
\item[(5)] $m'(s',u')=\infty$ for each $u'\in S'\setminus\{s',t'\}$, 
\item[(6)] $\phi(s)=w's't'(w')^{-1}$ and
\item[(7)] $\phi(t)=w't'(w')^{-1}$.
\end{enumerate}
\end{Theorem}

\begin{proof}
Suppose that $s\in S$ and 
$\phi(s)$ is not a reflection in $(W',S')$. 
Since $s^2=1$, $(\phi(s))^2=1$. 
By Lemma~\ref{lem3}, 
there exist $w',v'\in W'$ such that 
$\phi(s)=w'v'(w')^{-1}$ and 
$v'$ is the element of longest length in $W'_{S'(v')}$. 
Since $\phi(s)$ is not a reflection, $v'\not\in S'$, i.e., $|S'(v')|>1$. 
Hence $|S'(v')|=2$ because $\dim \Sigma(W',S')=2$. 
Let $S'(v')=\{s',t'\}$. 
Since $v'$ is the element of longest length in $W'_{S'(v')}=W'_{\{s',t'\}}$ 
and $v'$ is not a reflection, 
$m'(s',t')$ is even and 
$v'=(s't')^{m'}$, where $m'=m'(s',t')/2$.
Hence $\phi(s)=w'(s't')^{m'}(w')^{-1}$. 
By Lemma~\ref{lem5}, 
there exists a unique element $t\in S$ such that 
\begin{enumerate}
\item[(i)] $\phi(W_{\{s,t\}})=w'W'_{\{s',t'\}}(w')^{-1}$, 
\item[(ii)] $\phi(t)$ is a reflection in $(W',S')$, and 
\item[(iii)] $m(s,t)=m'(s',t')=2$. 
\end{enumerate}
Then $\phi(s)=w's't'(w')^{-1}$ by (iii).

Now $\phi(t)$ is a reflection by (ii) and 
\begin{align*}
\phi(t) &\in \phi(W_{\{s,t\}})=w'W'_{\{s',t'\}}(w')^{-1} \\
&=\{1,w's'(w')^{-1},w't'(w')^{-1},w's't'(w')^{-1}\}.
\end{align*}
Hence either $\phi(t)=w's'(w')^{-1}$ or $\phi(t)=w't'(w')^{-1}$. 
Here we may suppose that 
$$\phi(t)=w't'(w')^{-1}.$$

Finally we show that 
$m(s,u)=\infty$ for each $u\in S\setminus \{s,t\}$ and 
$m'(s',u')=\infty$ for each $u'\in S'\setminus \{s',t'\}$.

We suppose that there exists 
$u\in S \setminus \{s,t\}$ such that $m(s,u)<\infty$. 
By Lemma~\ref{lem4}, 
$\phi(W_{\{s,u\}})=x'W'_{\{a',b'\}}(x')^{-1}$
for some $x'\in W'$ and $a',b'\in S'$. 
Then 
$$ w's't'(w')^{-1}=\phi(s)\in 
\phi(W_{\{s,u\}})=x'W'_{\{a',b'\}}(x')^{-1}.$$
By Lemma~\ref{lem6}, 
$x'W'_{\{a',b'\}}(x')^{-1}=w'W'_{\{s',t'\}}(w')^{-1}$. 
Hence 
\begin{align*}
\phi(W_{\{s,u\}})&=x'W'_{\{a',b'\}}(x')^{-1} \\
&=w'W'_{\{s',t'\}}(w')^{-1} \\
&=\phi(W_{\{s,t\}}). 
\end{align*}
Thus $W_{\{s,u\}}=W_{\{s,t\}}$ and 
$\{s,u\}=\{s,t\}$ by Lemma~\ref{lem1}~(vi).
Hence $u=t$.
This contradicts the assumption $u\in S\setminus \{s,t\}$. 
Thus $m(s,u)=\infty$ for each $u\in S\setminus \{s,t\}$. 

We note that 
$$ \phi(st)=(w's't'(w')^{-1})(w't'(w')^{-1})=w's'(w')^{-1} $$
and 
$$ \phi^{-1}(s')=(\phi^{-1}(w'))^{-1}st \phi^{-1}(w').$$
By applying the above argument to $\phi^{-1}:W'\rightarrow W$, 
we can prove that 
$m'(s',u')=\infty$ for each $u'\in S'\setminus \{s',t'\}$. 
\end{proof}

We obtain the following theorem from Theorem~\ref{Thm1}.

\begin{Theorem}\label{Thm2}
Let $(W,S)$ and $(W,S')$ be Coxeter systems 
with two-dimensional Davis-Vinberg complexes. 
Then there exists $S''\subset W$ such that 
$(W,S'')$ is a Coxeter system which is isomorphic to $(W,S)$ 
and $R_{S'}=R_{S''}$.
\end{Theorem}

\begin{proof}
Let 
$$S_0=\{s\in S : \text{$s$ is not a reflection in $(W,S')$}\}
=\{s_1,\dots,s_n\}. $$
For each $i\in\{1,\dots,n\}$, 
there exists a unique element $t_i\in S\setminus S_0$ 
such that $m(s_i,t_i)=2$ by Theorem~\ref{Thm1}. 
Then $s_it_i$ is a reflection in $(W,S')$ by Theorem~\ref{Thm1}.
Let 
$$ S''=(S\setminus S_0)\cup\{s_1t_1,\dots,s_nt_n\}. $$
Then $(W,S'')$ is a Coxeter system which is isomorphic to $(W,S)$ 
by Lemma~\ref{lem7}.
Since $S''\subset R_{S'}$ by the constraction of $S''$, 
$R_{S''}=R_{S'}$ by Lemma~\ref{lem3-1}~(1).
\end{proof}

Theorem~\ref{Thm2} implies the following corollary.

\begin{Corollary}\label{Cor2}
For a Coxeter group $W$, 
if $(W,S)$ and $(W,S')$ are Coxeter systems 
with two-dimensional Davis-Vinberg complexes, then 
the Coxeter diagrams of $(W,S)$ and $(W,S')$ have 
the same number of vertices, the same number of edges 
and the same multiset of edge-labels.
\end{Corollary}

\begin{proof}
By Lemma~\ref{lem4}, 
the Coxeter diagrams of $(W,S)$ and $(W,S')$ have 
the same number of edges 
and the same multiset of edge-labels.
By Theorem~\ref{Thm2}, 
there exists $S''\subset W$ such that 
$(W,S'')$ is a Coxeter system which is isomorphic to $(W,S)$ 
and $R_{S'}=R_{S''}$.
Hence $|S|=|S''|=|S'|$ by Lemma~\ref{lem3-1}~(2).
\end{proof}

The main results of this paper 
have been announced at the Topology Symposium in Japan 
on July 20, 2003.

%

%
\end{document}